\documentclass[12pt]{amsart}
\usepackage{amsmath,amssymb,amsfonts,amscd,graphicx,xypic}
\usepackage{relsize}
\usepackage{exscale}

\usepackage{array,multirow}
\usepackage{color}

\usepackage{hyperref}
\hypersetup{
    colorlinks = true,
    linkcolor = {black},
   citecolor  = {black},
}

\newcommand{\rcong}{\mathbin{\rotatebox[origin=c]{-45}{$\cong$}}}

\usepackage[margin=1in]{geometry}
\usepackage{lmodern}
\usepackage{mathptmx}
\usepackage{amsthm}
\usepackage{amssymb}
\usepackage{tikz}
	\usetikzlibrary{decorations.pathreplacing}
	\usetikzlibrary{patterns}

\newcommand{\pr}{\prime}

\newcommand{\Z}{\mathbb{Z}}

\DeclareMathOperator{\wh}{Wh}
\DeclareMathOperator{\PP}{\, \text{P}}

\newcommand{\GoodMatrix}{
$
\begin{array}{c|cc|cc}
\multicolumn{5}{r}{ {{\scriptstyle{a_1}}}\;\;\;\;\; {{\scriptstyle{ b_1}}}     \;\;\;\;\, {{\scriptstyle{a_2}}}\;\;\;\;\; {{\scriptstyle{ b_2}}}}\\
  {{\scriptstyle{a_1}}} &  {0} &  {\pm1} &  {\star} &  {0} \\
 {{\scriptstyle{b_1}}} &  {0} &  {0} &  {0} &  {0} \\ 
 \hline
  {{\scriptstyle{a_2}}} &  {\star} &  {0} &  {0} &  {\pm1} \\
 {{\scriptstyle{b_2}}} &  {0} &  {0} &  {0} &  {0} \\ 
\end{array}
$
}

\usepackage{xcolor}

\usepackage{changes}

\newtheorem{theorem}[subsection]{Theorem}

\theoremstyle{definition}

\theoremstyle{remark}

\numberwithin{equation}{section}
\theoremstyle{plain}

\theoremstyle{definition}

\newcommand{\co}{\colon\thinspace}

\theoremstyle{definition}

\newtheorem*{defin*}{Definition}
\newtheorem{claim}{Claim}
\theoremstyle{plain}

\begin{document}

\title{Universal Surgery Problems with Trivial Lagrangian}

\author{Michael Freedman}
\address{\hskip-\parindent
  Michael Freedman\\
    Microsoft Research, Station Q, and Department of Mathematics\\
    University of California, Santa Barbara\\
    Santa Barbara, CA 93106}
\email{mfreedman@math.ucsb.edu}

\author{Vyacheslav Krushkal$^*$}
\thanks{$^*$Supported in part by NSF grant DMS-1612159.}
\address{\hskip-\parindent
  Slava Krushkal\\
    Department of Mathematics\\
    University of Virginia\\
    Charlottesville, VA 22904}
\email{krushkal@virginia.edu}
\begin{abstract} We study the effect of Nielsen moves and their geometric counterparts, handle slides, on good boundary links. A collection of links, universal for $4$-dimensional surgery, is shown to admit Seifert surfaces with trivial Lagrangian. 
They are good boundary links \cite{F01}, with Seifert matrices of a more general form than in known constructions of slice links. 
We show that a certain more restrictive condition on Seifert matrices is sufficient for proving the links are slice.  
We also give a correction of a Kirby calculus identity  in \cite{FK2}, useful for constructing surgery kernels associated to link-slice problems.
\end{abstract}

\maketitle

\section{Introduction}

This paper concerns the $4$-dimensional topological surgery conjecture, which is known to hold in the simply-connected case \cite{F} and more generally for a class of ``good'' fundamental groups. Its validity for arbitrary fundamental groups remains a central open problem. 
 
Universal surgery models may be formulated in terms of the free-slice problem for a collection of links in $S^3$. (A link is {\em freely slice} if the fundamental group of the slice complement in the $4$-ball is free, generated by meridians.) To describe the connection between surgery and link slicing problems in more detail, recall that
 a $k$-component link $L$ is a \emph{boundary link} if the components bound disjoint Seifert surfaces, or equivalently, there is a homomorphism  to the free group, ${\phi}\co \pi_1(S^3\backslash L) \twoheadrightarrow \text{Free}_k$, taking meridians to free generators. 
 $L$ is a \emph{good boundary link} \cite{F01} if $\text{ker}({\phi})$ is perfect for some $\phi$ as above. Good boundary links are known \cite{F01} to admit unobstructed surgery problems for constructing a slice complement. 
 
A stronger condition \cite{F93} is that in some symplectic basis $\{ a_1,\ldots, a_g, b_1, \ldots, b_g \}$ of simple closed curves for some choice of Seifert surfaces $S$,  the Seifert form is a direct sum of  blocks of the form
\begin{equation} \label{good}
\begin{array}{c|cc|}
\multicolumn{3}{r}{ {{\scriptstyle{a_i}}}\;\;\;\;\;\; {{\scriptstyle{ b_i}}}}\\
  {{\scriptstyle{a_i}}} &  {0} &  {\pm1} \\
 {{\scriptstyle{b_i}}} &  {0} &  {0}  \\ 
\multicolumn{3}{c}{}
\end{array}\;.
\end{equation}
For example, when $L$ has vanishing linking numbers, the Whitehead double of $L$ (with clasps of either sign) is of this type. 
Nielsen moves on ${\rm Free}_k$ have no effect on the good boundary property of a link. On the other hand,  condition (\ref{good}) and other more subtle, non-abelian, invariants of curves on $S$ are not ``gauge invariant" under band sums of Seifert surfaces, corresponding to Nielsen moves. 
In Theorem \ref{theorem} we show that geometric operations corresponding to Nielsen moves allow for a construction of links, universal for surgery, which admit Seifert surfaces with a ``trivial Lagrangian'', see details below and in section \ref{definitions section}. 

\begin{figure}[ht]
\includegraphics[width=14.3cm]{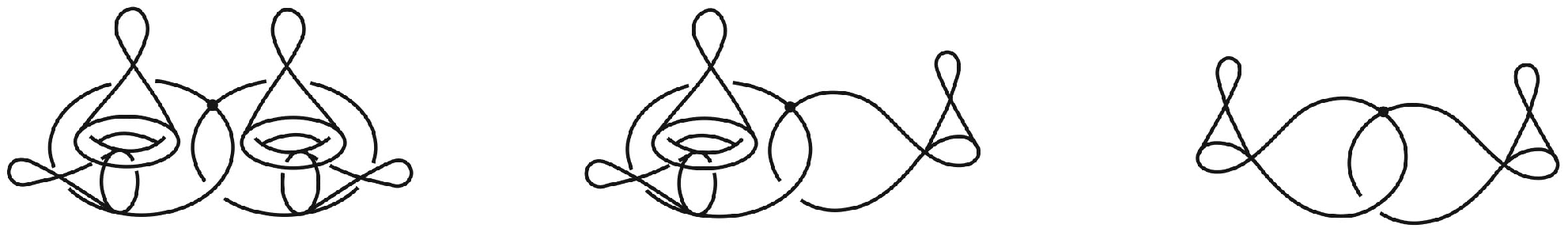}
{\small
\put(-410,-12){(1a) old universal problems}
\put(-262,-12){(1b) new universal problems}
\put(-80,-12){(1c) sufficient}
}
\caption{}
\label{Summary}
\end{figure}

Next we describe our results in more detail.
In \cite{FK1} a new class of universal surgery problems was produced where the surgery kernel is carried by an appropriately thickened 2-complex which appears closer to what is known to be sufficient to solve topological surgery problems. In icons, the progress is the middle picture in figure \ref{Summary}.

N.B. The icons ignore multiplicity of genus and double points which may make a real difference, e.g. the low multiplicity example pictured in Figure 1b is actually in the "sufficient" category, but its higher multiplicity cousins are not known to be.

In terms of the free link-slice problem,
the three stages may be summarized as the problem of slicing corresponding composite links, cf. \cite[Section 3]{FK1}:
{ $${ \text{\small (2a) } \,  \text{Wh}\circ \text{P}\circ  \text{Bing} \circ \text{P} \circ \text{Hopf}, \hspace{.2cm}  \text{\small (2b) }  \text{Wh}\circ \text{P}\circ \text{Bing} \circ \text{P} \circ  {\text{WhL}}, \hspace{.2cm}  \text{\small (2c) }\,  \text{Wh}\circ \text{P}\circ \text{Wh} \circ \text{P} \circ \text{Hopf},}$$}
where 
$$ \text{Hopf}\, =\, 
\vcenter{\hbox{
\includegraphics{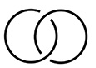}}},
 \; \;  {\text{WhL}}\, =\, 
\vcenter{\hbox{\includegraphics{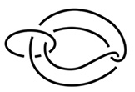}}}
, \; \text{Bing}\, =\, \vcenter{\hbox{\includegraphics{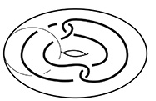}}}\, ,
\; \text{Wh}\, =\, \vcenter{\hbox{\includegraphics{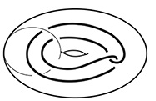}}}\, ,
$$
with either clasp in WhL and Wh. P denotes parallel copies, which allows for multiplicities.
General representatives of universal links of types 2a/2b are Whitehead doubles of homotopically essential links; as such they have usual genus one Seifert surfaces $S$ (cf. figures \ref{clasps}, \ref{SeifertSurfaces fig}). With this choice of Seifert surfaces, it is not difficult to see that any collection of simple closed curves, representing a Lagrangian subspace of $H_1(S;{\mathbb Z})$, forms a homotopically {\em essential} link.

Motivated, in part, by \cite{CKP} and the two questions raised in \cite[Section 7]{CKP}, we have found that all the universal surgery links of type 2b, after a suitable change of basis (corresponding to Nielsen moves on the free group) are ``good boundary links with a trivial Lagrangian" or ``Lagrangian-trivial" for short. This is notable in light of the main theorem of \cite{CKP} which shows that all links with the slightly stronger property "$\text{Lagrangian}$-trivial$^+$" are freely slice. All that stands in the way of a completely general topological surgery theorem, with only the high dimensional Wall obstruction, is the gap between Lagrangian-trivial and $\text{Lagrangian}$-trivial$^+$.

This gap could be real and the whole story (which we now suspect), or there may be some procedure consisting of Nielsen moves and clever choices of Seifert surfaces which allow a Lagrangian-trivial link to be promoted to a $\text{Lagrangian}$-trivial$^+$ link. Both possibilities will certainly be the subject of assiduous study.

In Section \ref{definitions section} we produce, in a simplified context, the Lagrangian-trivial link associated with cases 1b/2b. The simplification is that we actually exhibit the Nielsen moves from $\wh \circ \PP \circ \wh$ to a Lagrangian-trivial link. Skipping intermediate composition $\PP \circ \mathrm{Bing}$ in 2b is legitimate, as it corresponds to simply contracting the surface stages, see figure \ref{WhPWh fig} (compare with Figure 3.3 in \cite{FK1}). Free slicing the simplified link would imply a free slicing of case 2b. The 3-manifold that we actually analyze, the zero framed surgery ${\mathcal{S}}_{}^0(\wh \circ \PP \circ \rm{WhL})$, is the boundary of a 4-manifold with spine of the schematic form:
\begin{figure}[ht]
\includegraphics[height=2.cm, trim=7cm 0 0 0]{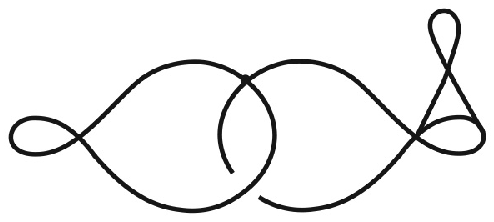}
\put(20,25){\text (with higher multiplicities allowed)}
\caption{}
\label{WhPWh fig}
\end{figure}

Section \ref{good*} states a condition on good boundary links, sufficient for solving the free-slice problem. It is of interest because the Seifert form is more general than that in (\ref{good}); specifically non-trivial linking of $a_i, a_j$ is allowed for $i\neq j$.
In Section \ref{correction section} we correct an error in our Kirby calculus, relevant to building slice complements \cite[Section 5.2]{FK2}, pointed out to us by the authors of \cite{CKP}.

\section{Main results} \label{definitions section}

Given a link $L$ with Seifert form, in some symplectic basis $\{ a_1,\ldots, a_g, b_1, \ldots, b_g \}$ of simple closed curves  on Seifert surfaces $S$,  a direct sum of  blocks of the form 
\begin{equation} \label{good1}
\begin{array}{c|cc|}
\multicolumn{3}{r}{ {{\scriptstyle{a_i}}}\;\;\;\;\;\; {{\scriptstyle{ b_i}}}}\\
  {{\scriptstyle{a_i}}} &  {0} &  {\pm1} \\
 {{\scriptstyle{b_i}}} &  {0} &  {0}  \\ 
\multicolumn{3}{c}{}
\end{array}\;,
\end{equation} there is a Lagrangian ($\frac{1}{2}$ dimensional) subspace of $H_1(S;\Z)$, for example the subspace spanned by $\{ b_1,\ldots, b_g\}$, on which linking and self-linking vanish. An even stronger condition \emph{Lagrangian-trivial} is that the Lagrangian subspace is spanned by (0-framed) disjoint simple closed curves $b_1, \dots, b_g$ which constitute a homotopically trivial link (i.e. all Milnor's $\overline{\mu}$-invariants with non-repeating indices vanish). Finally, the strongest condition considered here is \emph{$\text{Lagrangian}$-trivial$^+$} which requires that 
 the $2g$ long list of $(g+1)$-component links are each homotopically trivial. They are:
	\begin{align*}
		& \{b'_1 \cup \cdots \cup b'_g\} \cup a_i \hspace{1em} \forall i \hspace{1em} 1 \leq i \leq g, \text{ and }  \\
		& \{b'_1 \cup \cdots \cup b'_g\} \cup b_i \hspace{1em} \forall i \hspace{1em} 1 \leq i \leq g,
	\end{align*}
 {where $b'_i$ is a push-off  copy of $b_i$, having a trivial linking number with $a_i$.}
According to \cite{CKP},  {links that admit a good boundary basis $\{ a_i, b_i\}$ as above, satisfying the  $\text{Lagrangian}$-trivial$^+$ condition, are freely slice.} (Meridians in $\pi_1(S^3 \backslash L)$ map to free generators of $\pi_1$ of some topologically flat slice components.)  {We refer the reader to \cite{FK1} for a detailed discussion of the background material, including the notion of universal surgery problems, and Milnor's $\bar\mu$-invariants.} In this section we prove:

 \begin{theorem} \label{theorem} \sl The universal links of the form
$
\wh \circ \, \rm{P} \circ \rm{WhL}
$
are Lagrangian-trivial. 
\end{theorem}
$\rm P$ denotes any number of parallel copies taken on each side. In our discussion we treat in detail (and draw diagrams) only for the case of ${\rm P} = 2$ (on both components), which we write as ${\rm P}_{2,2}$. This is the first interesting case since this 4-component link is \emph{not} homotopically trivial so the natural Seifert surfaces do not exhibit Lagrangian-triviality. We achieve the Lagrangian-trivial property by a simple \emph{Nielsen move} on each side. We write the Whitehead link ${\rm WhL} = L_1 \cup L_2$. Then with parallels $\rm{P}_{2,2} \circ {\rm WhL} = L_1^0 \cup L_1^1 \cup L_2^0 \cup L_2^1$. The Nielsen moves are two handle slides which replace $\rm{P}_{2,2} \circ \wh$ with a 0-framed-surgery-equivalent link containing band-sums:
\[
(\rm{P} \circ {\rm WhL})^\pr : = L_1^0 \, \cup\,  L_1^0 \#_{\text{band}} L_1^1\,  \cup\,  L_2^0 \, \cup\,  L_2^0 \#_{\text{band}} L_2^1
\]

This is the case we draw, but if instead $\rm{P}_{k,j} \circ {\rm WhL} = L_1^0 \, \cup \cdots  \cup\,  L_1^k \, \cup\,  L_2^0 \, \cup \cdots \cup\,  L_2^j$, set
\[
(\rm{P} \circ {\rm WhL})^\pr = L_0^1 \, \cup\, L_1^0 \#_{\text{band}}L_1^1 \, \cup \cdots \cup \, L_1^0 \#_{\text{band}}L_1^k \# L_2^0 \, \cup \, L_2^0 \#_{\text{band}} L_2^1 \, \cup \cdots \cup \,  L_2^0 \#_{\text{band}} L_2^j.
\]

\begin{figure}[ht]
\includegraphics[height=2.3cm]{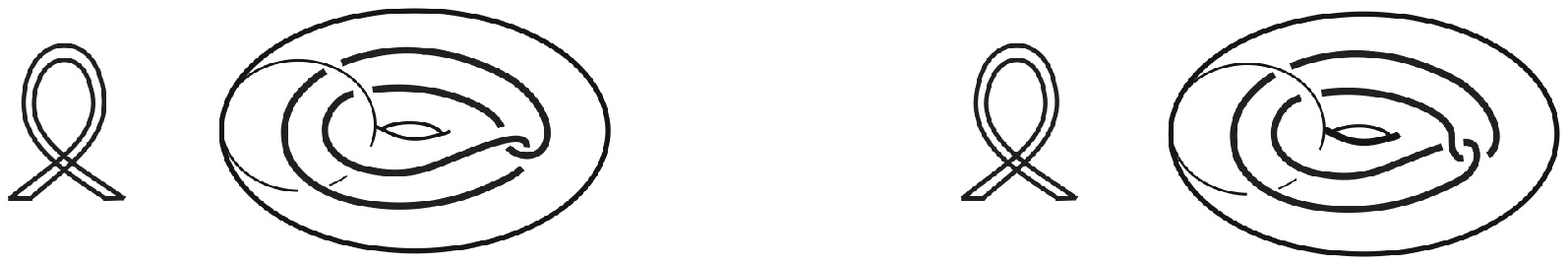}
\put(-390,25){$-$}
\put(-347,27){$\equiv$}
\put(-160,25){$+$}
\put(-118,27){$\equiv$}
{\small
\put(-256.5,16){$w_-$}
\put(-26.5,16){$w_+$}
}
\caption{}
\label{clasps}
\end{figure}

Before we start drawing pictures there is a small subtlety of double point signs to discuss. On the link level $+$ vs $-$ double points yield opposite clasps, figure \ref{clasps}.

In each case ($+$ or $-$) there are two pictures we draw (from Seifert's algorithm) for the Seifert surface $S^{\text{up}}_-$ or $S^{\text{down}}_-$ (resp. $S^{\text{up}}_+$ or $S^{\text{down}}_+$) bounding $w_-(w_+)$.
To make $S_\pm^\text{up (down)}$, plumb a $\pm$ twisted band above (below) with an annulus located in the $(x,y)$ plane.

\begin{figure}[ht]
\includegraphics[height=5.9cm]{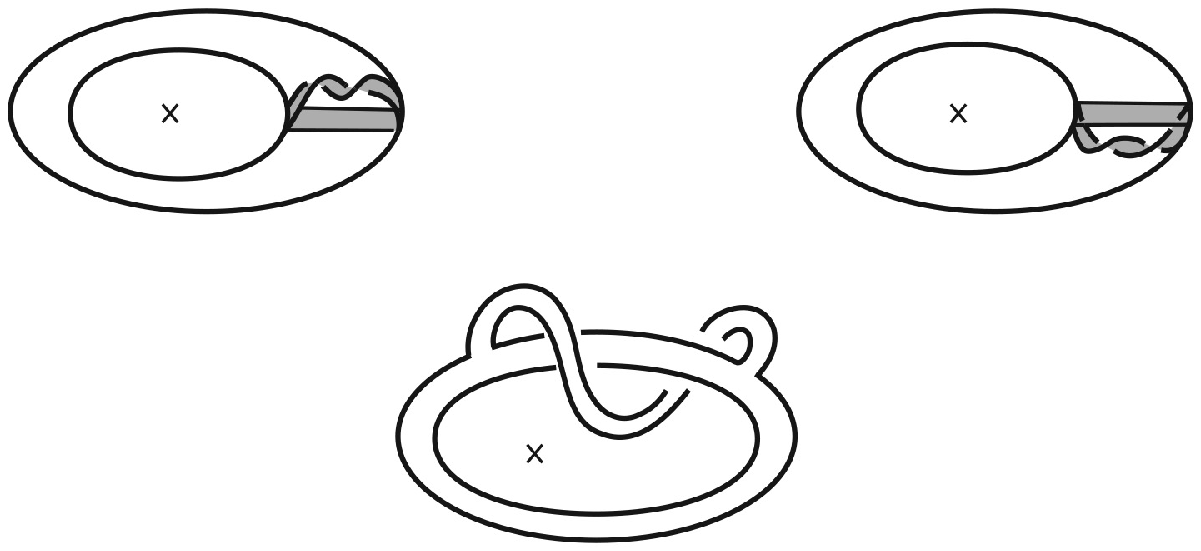}
\put(-367,130){$S_-^{\text{up}}$}
\put(-340,100){a}
\put(-12,100){b}
\put(-123,10){c}
\put(-3,135){$S_-^{\text{down}}$}
\caption{}
\label{SeifertSurfaces fig}
\end{figure}

$S_-^\text{up}$ and $S_-^\text{down}$ are pictured in figures \ref{SeifertSurfaces fig}a and \ref{SeifertSurfaces fig}b, respectively, with \ref{SeifertSurfaces fig}c an intermediate picture which is isotopic to both, showing $S_\pm^\text{up}$, $S_\pm^\text{down}$ are isotopic rel boundary; the apparent up/down choice for the band being merely a matter of basis choice in $H_1(S; \Z)$. The reason for distinguishing isotopic Seifert surfaces is that it will help us locate the correct band sum choices in the above expressions.

\begin{figure}[ht]
\includegraphics[height=6cm]{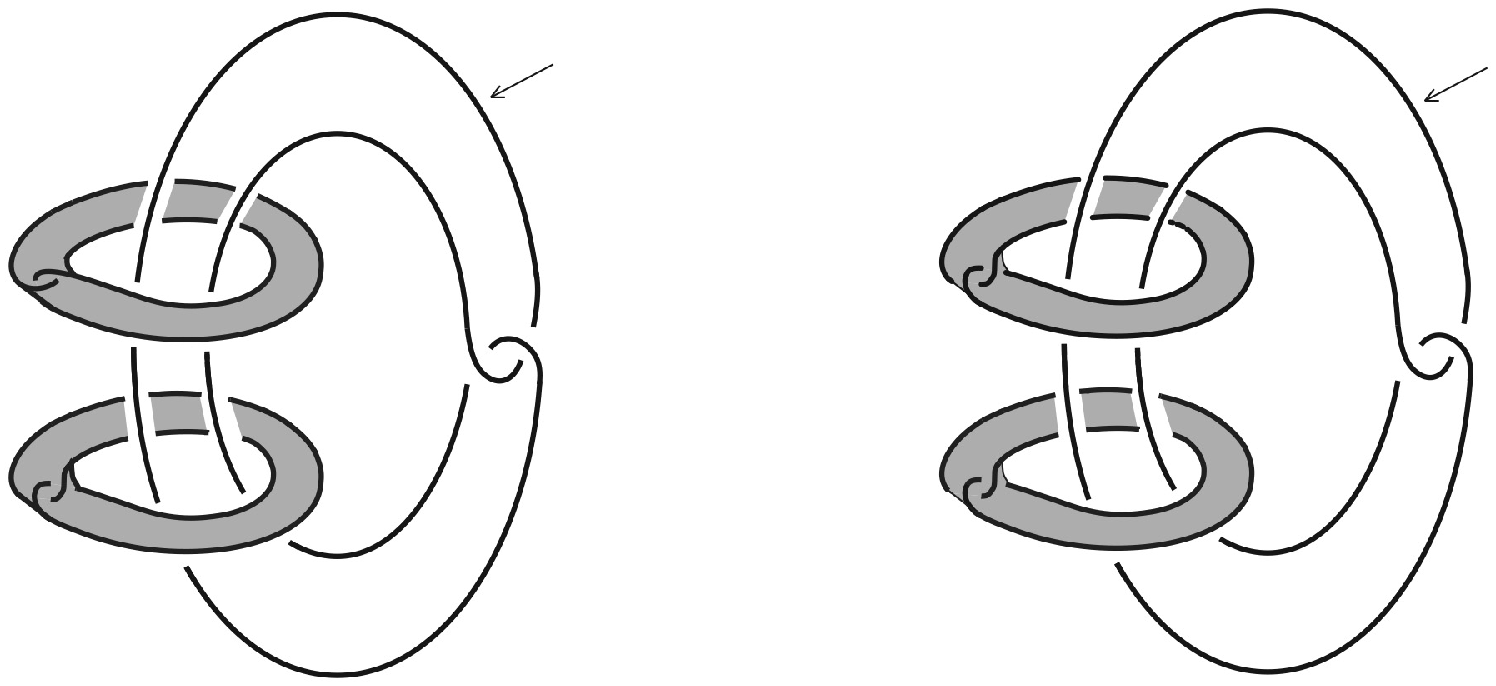}
{\small
\put(-373,107){up}
\put(-372,92){$+$}
\put(-383,53){down}
\put(-376,40){$-$}
\put(-348,25){$L^0_1$}
\put(-348,125){$L^1_1$}
\put(-128,25){$L^0_1$}
\put(-128,125){$L^1_1$}
\put(-166,107){down}
\put(-158,92){$-$}
\put(-166,53){down}
\put(-158,40){$-$}
\put(-238,15){$L_2$}
\put(-17,15){$L_2$}
\put(-1,140){${\rm Wh}\circ{\rm P}^{}_2$}
\put(-1,153){apply}
\put(-218,140){${\rm Wh}\circ{\rm P}^{}_2$}
\put(-218,153){apply}
}
\put(-205,75){or}
\put(-290,-15){Case 1}
\put(-70,-15){Case 2}
\caption{}
\label{WhP22}
\end{figure}

Now let us draw $\wh \circ \text{P}_{2,2} \circ {\rm WhL}$.
For clarity, we only draw the details of the link $L = \wh \circ {\rm P}_{2,2} \circ {\rm WhL}$ on the left side; we have also drawn Seifert surfaces for the final Whitehead doubles on the left and they should also be imagined on the right. A convention will help us get the Nielsen moves/band sums right. On each side we label one of the parallels $L_1^0$ (or $L_2^0$). In Figure \ref{WhP22}, $L_1^0$ is pictured as having a negative clasp, and a ``down" Seifert surface, but all that matters is that if $L_0^i$ has the opposite sign clasp from that of $L_0^0$, it should have the opposite kind (up vs. down) Seifert surface and if clasp sign is the same, then the same type of Seifert surface. These are cases 1 and 2 of Figure \ref{WhP22}.

\begin{figure}[ht]
\includegraphics[height=6.1cm]{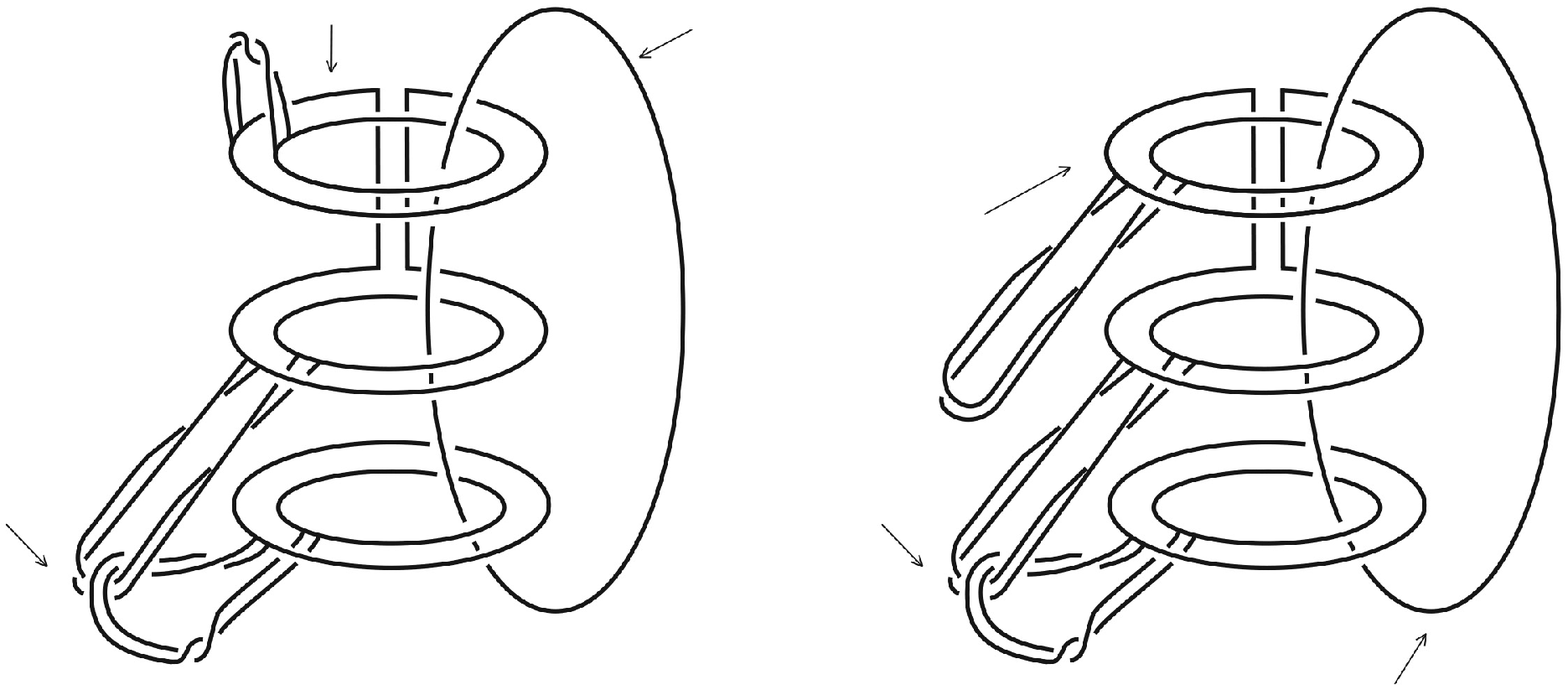}
{\small
\put(-70,0){apply}
\put(-218,165){apply}
\put(-321,174){$L_1^0 \#_{\text{band}}L_1^1$}
\put(-180,123){$L_1^0 \#_{\text{band}}L_1^1$}
}
{\scriptsize
\put(-80,-13){${\rm Wh}\circ{\rm P}^{}_2\circ {\rm Wh}$}
\put(-218,152){${\rm Wh}\circ{\rm P}^{}_2\circ {\rm Wh}$}
}
\put(-205,75){or}
\put(-403,43){$\circledast$}
\put(-181,44){$\circledast$}
\caption{Band-sum of Seifert surfaces corresponding to a Nielsen move. Lagrangian triviality$^+$ is foiled by linking marked with $\circledast$.}
\label{Seifert move}
\end{figure}

Now the passage to Figure \ref{Seifert move} post-composes the map $\pi_1(S^3-L) \twoheadrightarrow \mathrm{Free}$ with a Nielsen move, or more geometrically takes a parallel of $L_1^0$ and bands it to $L_1^1$, and on the level of Seifert surfaces bands a copy of the Seifert surface for $L_1^0$ to the Seifert surface for $L_1^1$.

\begin{figure}[ht]
\includegraphics[height=16cm]{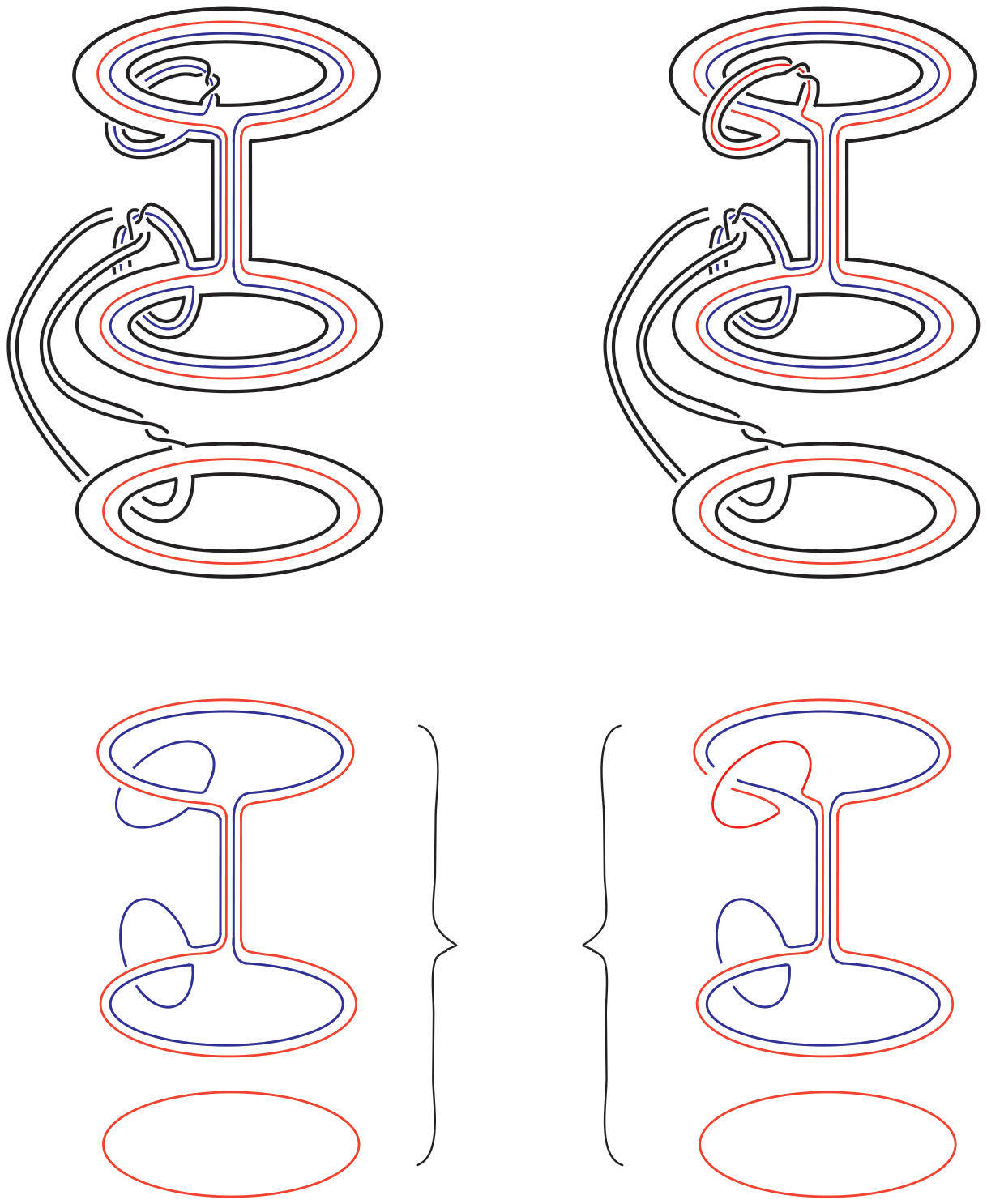}
\put(-182, 95){$\star$}
\caption{In both cases $\star$ shows a Lagrangian consisting of three $0$-framed components. The top two components are trivial in the solid torus, complementary of the vertical circle in Figure \ref{Seifert move}.}
\label{Seifert move1}
\end{figure}

In Figure \ref{Seifert move1} we have redrawn the genus $2$ Seifert surfaces to display more symmetry and indicated a $3$-component Lagrangian. The Lagrangian is for clarity displayed separately below in Figure \ref{Seifert move1}. Observe that two of the three components simply melt away, being $2$-component $0$-framed unlinks in the ambient solid torus, complementary to the vertical circle in Figure \ref{Seifert move}.

\begin{figure}[ht]
\includegraphics[height=3.5cm]{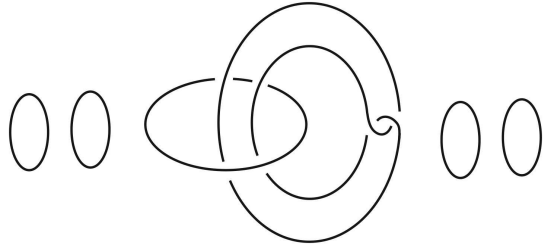}
\caption{The $6$-component Lagrangian.}
\label{AfterMoves}
\end{figure}

Now do precisely the same Nielsen move/band sum in a solid torus neighborhood of $L_2$ to obtain the same $3$-component Lagrangian on the right side. Combining the two sides, we see a $6$-component Lagrangian for $\wh \circ \text{P}_{2,2} \circ \rm{WhL}$ of the form in Figure \ref{AfterMoves}.

The sign of the indicated clasp is the sign of the clasp of WhL in $\wh \circ \text{P}_{2,2} \circ \rm{WhL}$. The outer Whitehead doublings have no influence on the ultimate Lagrangian link.  The link in Figure \ref{AfterMoves} is easily seen to be homotopically trivial, proving that $\wh \circ {\rm P}_{2,2} \circ \wh$ is \emph{Lagrangian-trivial}. The linking labeled $\circledast$
in Figure \ref{Seifert move} implies that these Seifert surfaces do \emph{not} exhibit Lagrangian triviality$^+$.

As noted above, the general case $\wh \circ {\rm P}_{k, j} \circ \wh$ is analogous. In this case Figure \ref{AfterMoves} becomes a Whitehead link with $2(k+j)$ additional unlinked unknotted components.
\qed

\subsection{Slice links} \label{good*}
\begin{figure}[ht]
\includegraphics[height=5.5cm, trim=7.5cm 0 0 0]{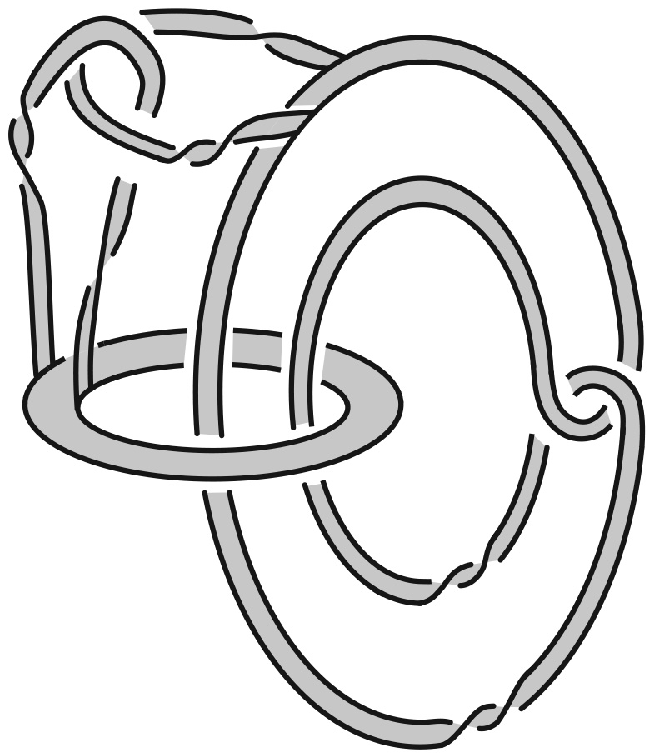} 
 \put(50,80){\GoodMatrix}
 \put(-170,70){$L$}
 \caption{An example of a link and the corresponding Seifert matrix. Here $b_1, b_2$ are the cores of the untwisted bands, thickening the Whitehead link. $a_1, a_2$ are the $(1,1)$ curves in the indicated genus $1$ Seifert surfaces. For the pictured link, the entries labeled $\star$ in the matrix are $\pm 1$.}
\label{WhPrime}
\end{figure}

It is shown in \cite{CKP} that links admitting a good boundary basis (the Seifert matrix is a direct sum of the form (\ref{good1})), satisfying the  $\text{Lagrangian}$-trivial$^+$ condition, are freely slice. The trivial Lagrangian, constructed in the proof of Theorem \ref{theorem}, is not trivial$^+$   because of the linking indicated by the symbol $\circledast$ in figure \ref{Seifert move}. In particular, the Seifert matrix has additional non-zero entries, corresponding to linking of $a_i, b_j$ for $i\neq j$. In this section we show that the result of \cite{CKP} extends to the setting of non-trivial linking numbers ${\rm lk}(a_i, a_j)$, $i\neq j$. An example of this type is shown in figure \ref{WhPrime}. The Seifert surfaces are obtained by plumbing untwisted bands, thickening the Whitehead link, and $\pm 1$ twisted bands which link as indicated in the figure. This link $L$  is akin to the Whitehead double of the Whitehead link, except that the twisted bands link. $L$ is Lagrangian-trivial$^+$ since all four links used in the verification of this condition are homotopy-trivial, being the Whitehead link with a parallel copy of one of its components.

Specifically, we show: {\sl Suppose the components of a link $L$ have disjoint Seifert surfaces with a symplectic basis
$\{a_1, \ldots, a_g, b_1,\ldots, b_g\}$ of simple closed curves, which is trivial$^+$ as defined in the beginning of section \ref{definitions section}. Suppose also that the Seifert matrix in this basis has diagonal blocks of the form (\ref{good1}), and the off-diagonal entries of the form ${\rm lk}(a_i, b_j)$, $i\neq j$ and ${\rm lk}(b_i, b_j)$, $i\neq j$ are all zero. Then $L$ is freely slice.}

To prove this statement, we summarize the argument in the case where the Seifert matrix is a direct sum of the blocks (\ref{good1}) (see \cite{CKP} for details), and indicate how to complete the proof when the linking numbers ${\rm lk}(a_i,a_j),\,  i\neq j$ are not necessarily zero. Consider the $4$-manifold $W$, obtained by attaching round $1$-handles to $D^4$ along $b^+_i$ and $b^-_i$, and zero-framed $2$-handles along $a_i$, for each $i$. Here $b^+_i$, $b^-_i$ are $+, -$ push-offs of $b_i$ in the normal direction to the Seifert surface. The boundary of $W$ is diffeomorphic to ${\mathcal S}^0(L)$, the zero-framed surgery on $L$ \cite{F93}. The surgery kernel is represented by hyperbolic pairs (torus, sphere), where the tori are formed by cores of the round $1$-handles and annuli bounded by $b^+_i$ and $b^-_i$ in $S^3$, and the spheres $A_i$ are formed by cores of the zero-framed $2$-handles, capped off by null-homotopies of the curves $\{ a_i\} $ in $D^4$. As in \cite{CKP} (building on earlier work \cite{FT}), the trivial$^+$ condition is used to construct a collection of singular disks in $D^4$, bounded by $\{ a_i, b^+_i, b^-_i\}$, such that the entire collection of $2$-spheres $K:=\cup_{i=1}^g (A_i\cup B_i)$ is ${\pi}_1$-null. Here $\{ B_i\}$ are obtained by capping off the cores of the round handles by disks bounded by $b^+_i, b^-_i$. When linking numbers ${\rm lk}(a_i,a_j),\,  i\neq j$ are non-zero, the spheres $A_i, A_j$ have non-trivial algebraic intersections. Use Norman's moves on these spheres (tubing $A_i$ into a parallel copy of $B_j$)  to achieve trivial algebraic intersections of $A_i, A_j$. The $2$-complex $K$ is ${\pi}_1$-null, so the Norman move, taking place in a neighborhood of $K$, preserves the ${\pi}_1$-null condition.
Now the proof is completed by \cite[Theorem 6.1]{FQ}, producing embedded spheres to complete surgery. \qed

\section{Correction} \label{correction section}

A useful tool for constructing surgery kernels associated to link-slice problems is an identity describing zero framed surgery $\mathcal{S}^0(L)$ on a \emph{good boundary} link $L$ as the result of \emph{plumbed Lagrangian surgeries} and additional \emph{twist surgeries}. In \cite[Lemma 5.8]{FK2} we published a version of this identity which erroneously lacked the $\pm 1$ twist surgeries. We would like to thank the authors of \cite{CKP} for calling this to our attention, and here we publish the correct identity. First the genus $2$ case.

\begin{figure}[ht]
\includegraphics[height=4.8cm]{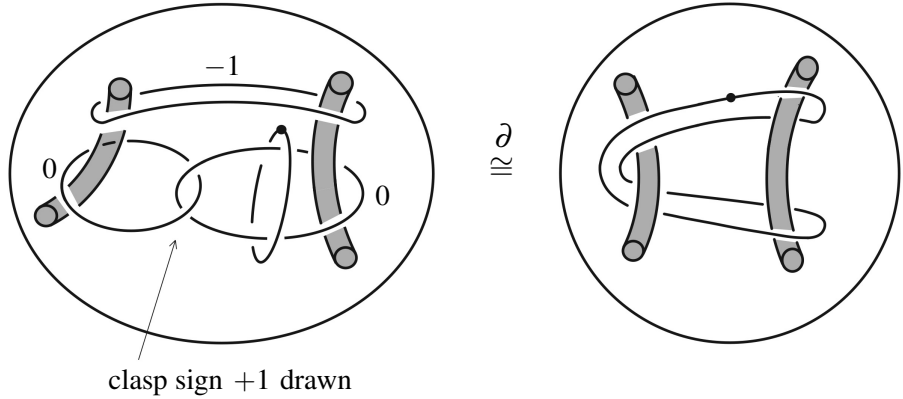} 
{\small
\put(-265,109){$-1$}
\put(-200,60){$0$}
\put(-326,70){$0$}
\put(-300,-10){clasp sign $+1$ drawn}
}
\put(-157,70){$\cong$}
\put(-155,82){$\partial$}
 \caption{Shaded tubes are removed from the balls, creating genus. An analogous result holds for clasp sign $-1$.}
\label{Calc1}
\end{figure}

\begin{claim} \label{Claim1}
	The two sides of Figure \ref{Calc1} are surgery diagrams which represent diffeomorphic 3-manifolds, the diffeomorphism being the identity on the boundary of the solid genus 2 handelbody containing the two diagrams.
\end{claim}

\begin{figure}[ht]
\includegraphics[height=10.5cm]{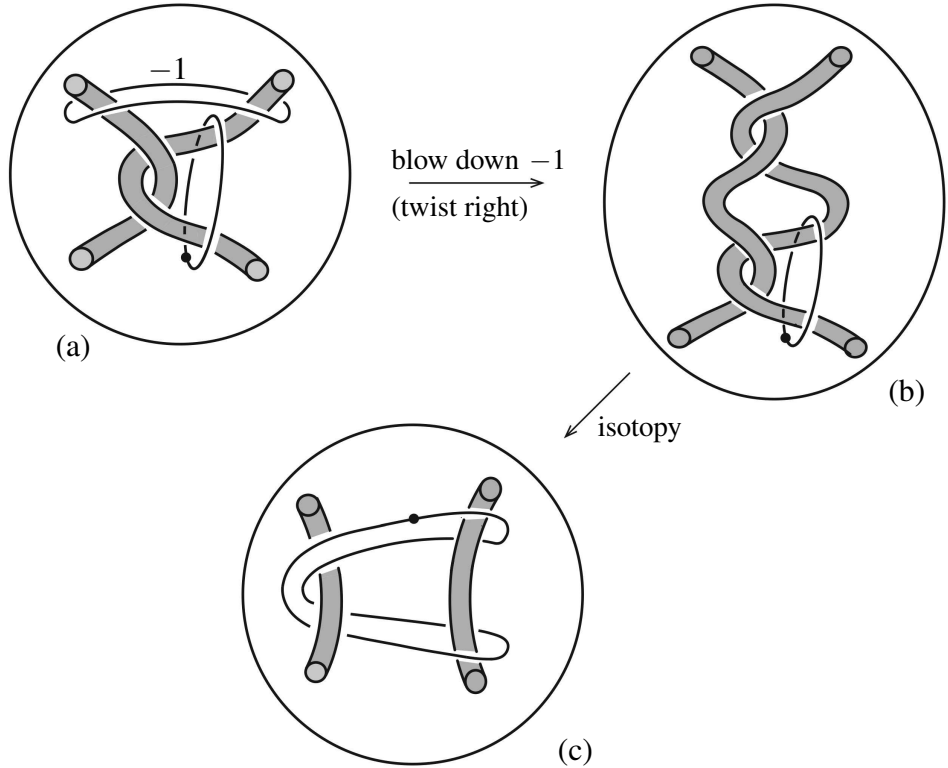} 
\put(-303,265){$-1$}
\put(-337,160){(a)}
\put(-22,143){(b)}
\put(-147,7){(c)}
{\small
\put(-210,231){blow down $-1$}
\put(-210,213){(twist right)}
\put(-132,130){isotopy}
}
 \caption{Proof of Claim \ref{Claim1}.}
\label{Calc2}
\end{figure}

\begin{proof}
	Canceling the plumbed pair yields Figure \ref{Calc2} (a); the rest of the calculation is given in Figure \ref{Calc2}.
\end{proof}

\begin{figure}[ht]
\includegraphics[width=14cm]{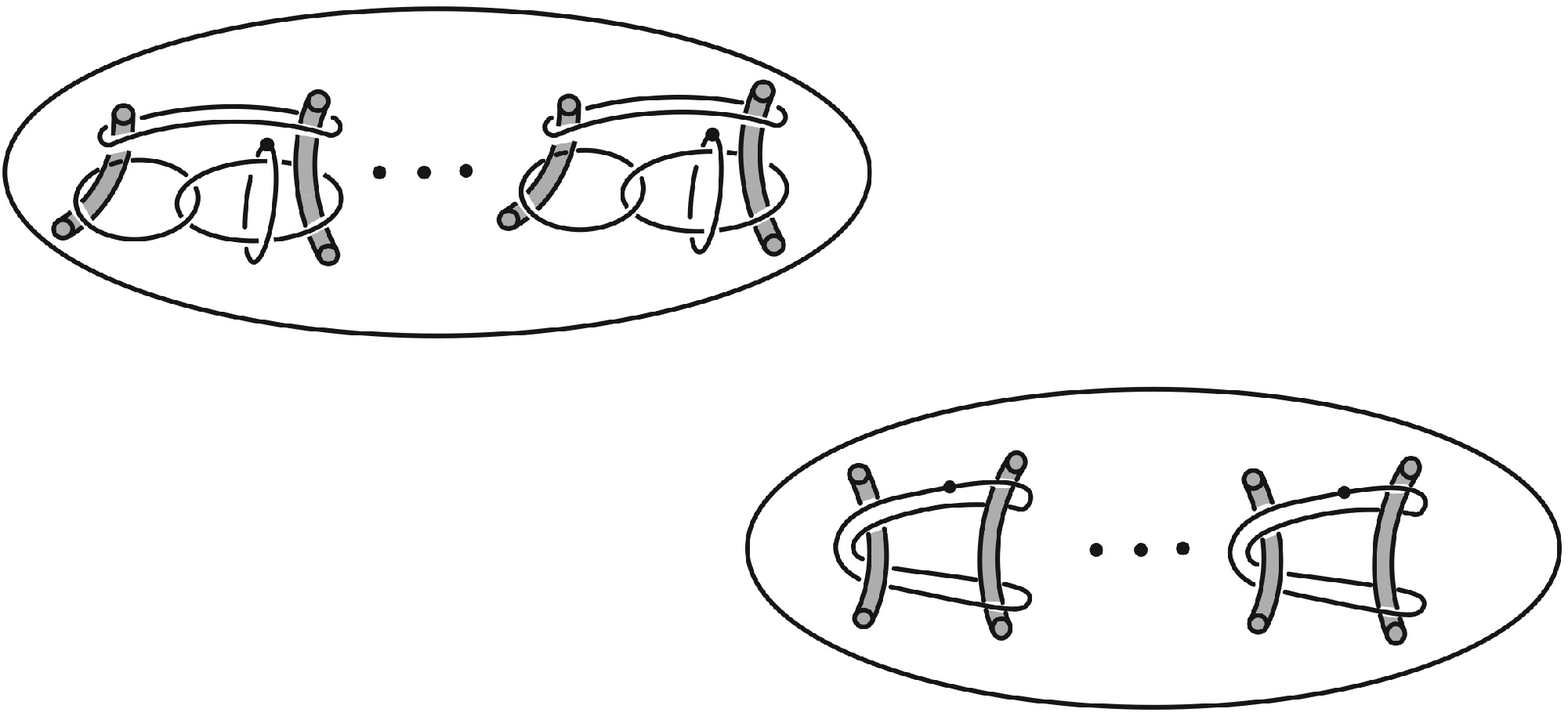} 
{\scriptsize
\put(-350,161){$-1$}
\put(-240,163){$-1$}
\put(-388,135){$0$}
\put(-310,129){$0$}
\put(-260,117){$0$}
\put(-198,134){$0$}
}
{\large
\put(-207,85){$\rcong$}}
\put(-194,94){$\partial$}
 \caption{}
\label{Calc3}
\end{figure}

\begin{claim}
	By the same argument in genus $= 2n$ there is an identity in figure \ref{Calc3}.
\end{claim}

\begin{figure}[ht]
\includegraphics[height=7.2cm]{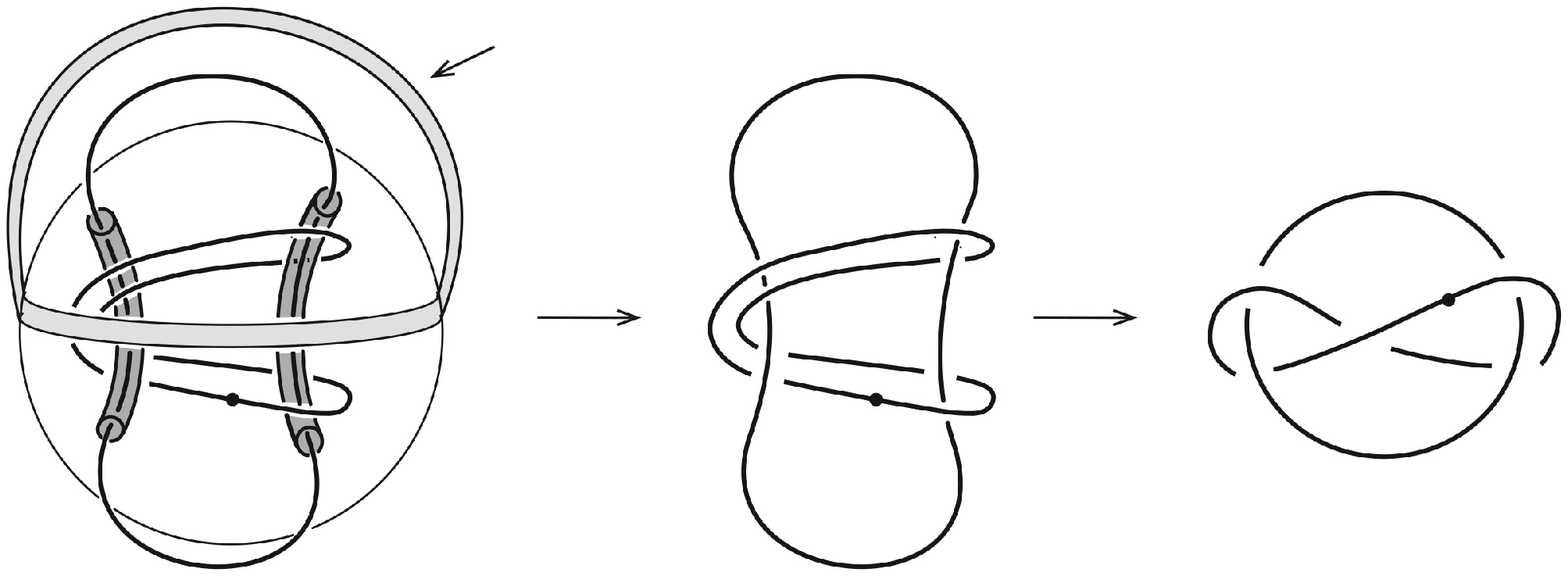} 
\put(-170,10){(b)}
\put(-100,20){(c) Whitehead link}
\put(-333,10){(a)}
 \put(-290,151){$2$-handle}
  \caption{}
\label{Calc4}
\end{figure}

Recall in the proof of  \cite[Lemma 5.8]{FK2} we had localized the calculation to genus $2$-handlebody pictured in Figure \ref{Calc2}a (but with the $\pm 1$ framed simple closed curve missing). The dotted simple closed curve in Figure \ref{Calc2}c bounds an obvious genus one Seifert surface in the handlebody. Note that the dotted curve becomes the Whitehead double bounding this genus one Seifert surface within a solid torus obtained from the original handlebody by attaching a $3$-dimensional $2$-handle to the equator of Figure \ref{Calc2}c. This $2$-handle, if present, is equivalent to the two unlinked Lagrangian curves $x_i$ and $y_i$ (see \cite{FK2}) being parallel, in which case the generalized double is, in fact, an ordinary Whitehead double. We illustrate this in Figure \ref{Calc4}.

Note that in this case of an honest Whitehead double the $\pm 1$ framed simple closed curve ``slips off" and does not affect the calculation. This is why it was overlooked in \cite{FK2}.

Our correction applied to Lemmas 5.8 and 5.10 of \cite{FK2} says that the $4$-manifolds which are constructed there are not (unobstructed) surgery problems for the original generalized Whitehead double link $\text{GWD} \subset S^3$, but rather for a related link in some integral homology sphere $\Sigma^3$ which is the result of the (overlooked) $\pm 1$ framed simple closed curves.

\end{document}